\documentclass[12pt]{amsart}
\usepackage{mathtools}
\begin{document}
\numberwithin{equation}{section}

\def\1#1{\overline{#1}}
\def\2#1{\widetilde{#1}}
\def\3#1{\widehat{#1}}
\def\4#1{\mathbb{#1}}
\def\5#1{\frak{#1}}
\def\6#1{{\mathcal{#1}}}

\newcommand{\UH}{\mathbb{H}}
\newcommand{\de}{\partial}
\newcommand{\R}{\mathbb R}
\newcommand{\Ha}{\mathbb H}
\newcommand{\al}{\alpha}
\newcommand{\tr}{\widetilde{\rho}}
\newcommand{\tz}{\widetilde{\zeta}}
\newcommand{\tk}{\widetilde{C}}
\newcommand{\tv}{\widetilde{\varphi}}
\newcommand{\hv}{\hat{\varphi}}
\newcommand{\tu}{\tilde{u}}
\newcommand{\tF}{\tilde{F}}
\newcommand{\debar}{\overline{\de}}
\newcommand{\Z}{\mathbb Z}
\newcommand{\C}{\mathbb C}
\newcommand{\Po}{\mathbb P}
\newcommand{\zbar}{\overline{z}}
\newcommand{\G}{\mathcal{G}}
\newcommand{\So}{\mathcal{S}}
\newcommand{\Ko}{\mathcal{K}}
\newcommand{\U}{\mathcal{U}}
\newcommand{\B}{\mathbb B}
\newcommand{\oB}{\overline{\mathbb B}}
\newcommand{\Cur}{\mathcal D}
\newcommand{\Dis}{\mathcal Dis}
\newcommand{\Levi}{\mathcal L}
\newcommand{\SP}{\mathcal SP}
\newcommand{\Sp}{\mathcal Q}
\newcommand{\A}{\mathcal O^{k+\alpha}(\overline{\mathbb D},\C^n)}
\newcommand{\CA}{\mathcal C^{k+\alpha}(\de{\mathbb D},\C^n)}
\newcommand{\Ma}{\mathcal M}
\newcommand{\Ac}{\mathcal O^{k+\alpha}(\overline{\mathbb D},\C^{n}\times\C^{n-1})}
\newcommand{\Acc}{\mathcal O^{k-1+\alpha}(\overline{\mathbb D},\C)}
\newcommand{\Acr}{\mathcal O^{k+\alpha}(\overline{\mathbb D},\R^{n})}
\newcommand{\Co}{\mathcal C}
\newcommand{\Hol}{{\sf Hol}}
\newcommand{\Aut}{{\sf Aut}(\mathbb D)}
\newcommand{\D}{\mathbb D}
\newcommand{\oD}{\overline{\mathbb D}}
\newcommand{\oX}{\overline{X}}
\newcommand{\loc}{L^1_{\rm{loc}}}
\newcommand{\la}{\langle}
\newcommand{\ra}{\rangle}
\newcommand{\thh}{\tilde{h}}
\newcommand{\N}{\mathbb N}
\newcommand{\kd}{\kappa_D}
\newcommand{\Hr}{\mathbb H}
\newcommand{\ps}{{\sf Psh}}
\newcommand{\Hess}{{\sf Hess}}
\newcommand{\subh}{{\sf subh}}
\newcommand{\harm}{{\sf harm}}
\newcommand{\ph}{{\sf Ph}}
\newcommand{\tl}{\tilde{\lambda}}
\newcommand{\gdot}{\stackrel{\cdot}{g}}
\newcommand{\gddot}{\stackrel{\cdot\cdot}{g}}
\newcommand{\fdot}{\stackrel{\cdot}{f}}
\newcommand{\fddot}{\stackrel{\cdot\cdot}{f}}

\def\Re{{\sf Re}\,}
\def\Im{{\sf Im}\,}

\newcommand{\Real}{\mathbb{R}}
\newcommand{\Natural}{\mathbb{N}}
\newcommand{\Complex}{\mathbb{C}}
\newcommand{\ComplexE}{\overline{\mathbb{C}}}
\newcommand{\Int}{\mathbb{Z}}
\newcommand{\UD}{\mathbb{D}}
\newcommand{\clS}{\mathcal{S}}
\newcommand{\gtz}{\ge0}
\newcommand{\gt}{\ge}
\newcommand{\lt}{\le}
\newcommand{\fami}[1]{(#1_{s,t})}
\newcommand{\famc}[1]{(#1_t)}
\newcommand{\ts}{t\gt s\gtz}
\newcommand{\classCC}{\tilde{\mathcal C}}
\newcommand{\classS}{\mathcal S}

\newcommand{\Step}[2]{\begin{itemize}\item[{\bf Step~#1.}]{\it #2}\end{itemize}}
\newcommand{\step}[2]{\begin{itemize}\item[{\it Step~#1.}]{\it #2}\end{itemize}}
\newcommand{\proofbox}{\hfill$\Box$}
\newcommand{\claim}[2]{\begin{itemize}\item[{\it Claim~#1.}]{\it #2}\end{itemize}}

\newcommand{\mcite}[1]{\csname b@#1\endcsname}
\newcommand{\UC}{\mathbb{T}}

\newcommand{\Moeb}{\mathrm{M\ddot ob}}

\newcommand{\dAlg}{{\mathcal A}(\UD)}
\newcommand{\diam}{\mathrm{diam}}

\def\mydot{\stackrel{\text{\Huge.}}}
\def\Slope{\mathop{\mathsf{Slope}}\nolimits}

\theoremstyle{theorem}
\newtheorem {result} {Theorem}
\setcounter {result} {64}
 \renewcommand{\theresult}{\char\arabic{result}}



\newcommand{\Spec}{\Lambda^d}
\newcommand{\SpecR}{\Lambda^d_R}
\newcommand{\Prend}{\mathrm P}




\def\cn{{\C^n}}
\def\cnn{{\C^{n'}}}
\def\ocn{\2{\C^n}}
\def\ocnn{\2{\C^{n'}}}
\def\je{{\6J}}
\def\jep{{\6J}_{p,p'}}
\def\th{\tilde{h}}


\def\dist{{\mathsf{dist}}}
\def\const{{\rm const}}
\def\rk{{\rm rank\,}}
\def\id{{\sf id}}
\def\aut{{\sf aut}}
\def\Aut{{\sf Aut}}
\def\CR{{\rm CR}}
\def\GL{{\sf GL}}
\def\Re{{\sf Re}\,}
\def\Im{{\sf Im}\,}
\def\U{{\sf U}}

\def\la{\langle}
\def\ra{\rangle}

\newcommand{\sgn}{\mathop{\mathrm{sgn}}}

\emergencystretch15pt \frenchspacing

\newtheorem{theorem}{Theorem}[section]
\newtheorem{lemma}[theorem]{Lemma}
\newtheorem{proposition}[theorem]{Proposition}
\newtheorem{corollary}[theorem]{Corollary}
\newtheorem{conjecture}{Conjecture}

\theoremstyle{definition}
\newtheorem{diff}{Difference}
\newtheorem{definition}[theorem]{Definition}
\newtheorem{notation}[theorem]{Notation}
\newtheorem{example}[theorem]{Example}

\theoremstyle{remark}
\newtheorem{remark}[theorem]{Remark}
\numberwithin{equation}{section}

\newcommand{\Maponto}
{\xrightarrow{\hbox{\lower.2ex\hbox{$\scriptstyle \smash{\mathsf{onto}}$}}\,}}
\newcommand{\Mapinto}
{\xrightarrow{\hbox{\lower.2ex\hbox{$\scriptstyle \smash{\mathsf{into}}$}}\,}}
\newcommand{\anglim}{\angle\lim}

\newenvironment{mylist}{\begin{list}{}%
{\labelwidth=2em\leftmargin=\labelwidth\itemsep=.4ex plus.1ex
minus.1ex\topsep=.7ex plus.3ex
minus.2ex}%
\let\itm=\item\def\item[##1]{\itm[{\rm ##1}]}}{\end{list}}

\title[Slope problem for semigroups]{Slope problem for trajectories of holomorphic semigroups in the unit disc}

\author[M. D. Contreras]{Manuel D. Contreras $^\dag$}

\author[S. D\'{\i}az-Madrigal]{Santiago D\'{\i}az-Madrigal $^\dag$}
\address{Camino de los Descubrimientos, s/n\\
Departamento de Matem\'{a}tica Aplicada II and IMUS\\
Universidad de Sevilla\\
Sevilla, 41092\\
Spain.}\email{contreras@us.es} \email{madrigal@us.es}

\author[P. Gumenyuk]{Pavel Gumenyuk $^\ddag$}
\address{
Dipartimento di Matematica \\
Universit\`a di Roma "Tor Vergata" \\
Via della Ricerca Scientifica 1 \\
00133 Roma, ITALIA.}
\email{gumenyuk@mat.uniroma2.it}

\date{\today}


\thanks{$^\dag$ Partially supported by the \textit{Ministerio
de Econom\'{\i}a y Competitividad} and the European Union (FEDER), project MTM2012-37436-C02-01, and  by \textit{La Consejer\'{\i}a de Educaci\'{o}n y Ciencia de la Junta de Andaluc\'{\i}a}.}

\thanks{$^\ddag$ Partially supported by the FIRB grant Futuro in Ricerca ``Geometria Differenziale Complessa e Dinamica Olomorfa'' n. RBFR08B2HY}

\begin{abstract}
It has been an open problem for about ten years whether every trajectory of a parabolic one-parameter semigroup in the unit disk tends to the Denjoy\,--\,Wolff point with a definite (and common for all trajectories) slope. In this paper, we give the negative answer to this question.
\end{abstract}

\maketitle

\def\Gen{\mathcal G}
\def\Gena{\Gen_{\mathrm{Aut}}}
\def\RH{\mathbb C_+}
\def\mes{\mathop{\mathrm{mes}}}

\newcommand{\REM}[1]{\relax}
\newcommand{\blogger}[2]{\vskip3mm\noindent{\large\underline{\textbf{\textsc{\uppercase{#1}}, \textit{#2}}}}}

\section{Introduction}
A \textsl{one-parameter semigroup of holomorphic functions} in the unit disc~$\UD:=\{z\colon|z|<1\}$ is a continuous
homomorphism $t\mapsto\varphi_t$ between the semigroup of non-negative real numbers endowed with the Euclidean topology and the semigroup~$\Hol(\UD,\UD)$ of all holomorphic self-maps of the unit disc endowed with the topology of uniform convergence on compacta.

The study of this kind of one-parameter semigroups has a long history. They have appeared in many different areas of Analysis such as operator theory, geometric function theory, optimization and branching stochastic processes, see, \textit{e.g.}, \cite{Goryainov::survey}. Moreover, one-parameter semigroups of holomorphic functions can be regarded as a special but very important case (namely, the autonomous case from the dynamical point of view) of a more general notion playing a fundamental role in Loewner theory, a current and very active branch of research (see, \textit{e.g.},  the recent survey \cite{ABCD}). We refer the reader to the monographs \cite{Shoikhet-book} and \cite{ES-book} for the basic definitions and results of the theory of one-parameter semigroups.

In the last decades, starting with the seminal paper by Berkson and Porta~\cite{BP}, the study of the dynamical aspects of the theory has attracted considerable interest. Namely, for any one-parameter semigroup $(\varphi _{t})\subset\Hol(\UD,\UD)$, there exists a \textit{unique}
analytic function $G:\Bbb{D\rightarrow C}$ (the so-called vector field or \textsl{infinitesimal generator} of the semigroup) such that, for each $z\in \Bbb{D}$, the positive trajectory associated with the Cauchy problem
\begin{equation*}
\left\{
\begin{array}{l}
\mydot{w}=G(w), \\
w(0)=z,
\end{array}
\right.
\end{equation*}
is exactly $[0,+\infty)\ni t\mapsto \varphi_t(z)$. This fact gives rise, in a natural way, to a number of dynamical questions concerning the complete (not only the positive) trajectories of the vector field~$G$ and thus leading to dynamical questions for the corresponding semigroup~$(\varphi_t)$. This dynamical approach includes topics like the asymptotic behavior of the trajectories, their $\omega $-limits, the analysis of the boundary fixed points, the multipliers of these points or the slopes of the trajectories arriving to a boundary fixed point.
 It is worth mentioning that this dynamical point of view (in particular, analysis of slopes) has also been extensively treated for holomorphic discrete iteration in the unit disk (see, \textit{e.g}, the papers by Poggi-Corradini \cite{Poggi-Corradini} and Pommerenke \cite{Pommerenke}).

We recall that a point $\xi$ in the Riemann Sphere $ \Bbb{C}_{\infty }$ is called an \textsl{$\omega $-limit point} of a curve $\Gamma
:(s_{1},s_{2})\rightarrow \Bbb{C}$, $-\infty < s_{1}<s_{2}\leq +\infty$,
if there exists a strictly increasing sequence $(t_{n})\subset (s_{1},s_{2})$
convergent to $s_{2}$ such that $\Gamma \left( t_{n}\right) \rightarrow \xi$. The set of all $\omega $-limit points of $\Gamma $ is called the \textsl{$\omega$-limit set} and denoted by $\omega \left( \Gamma \right)$.

Concerning the asymptotic behavior, the most fundamental result is the continuous version of the Denjoy\,--\,Wolff theorem, see, \textit{e.g.}, \cite[Proposition 4.1.2]{Shoikhet-book}, stating basically that for any one-parameter semigroup $(\varphi_t)\subset\Hol(\UD,\UD)$ such that $\varphi_{t_0}$ has no fixed point in~$\UD$ for some~${t_0>0}$, there exists a (unique) point $\tau\in\partial\D$, called the \textsl{Denjoy\,--\,Wolff point} of~$(\varphi_t)$, such that
\[
\lim_{t\rightarrow+\infty}\varphi_t(z)=\tau\quad \textrm{ for every } z\in\D.
\]
Moreover,  see, \textit{e.g.}, \cite[Theorem 2.12]{Shoikhet-book}, $\lim_{r\rightarrow1^{-}}\varphi_t^\prime(r\tau)=\exp(-\lambda t)$ for some $\lambda\ge0$ and all~$t\ge0$. The semigroup~$(\varphi_t)$ is said to be \textsl{hyperbolic} if~$\lambda>0$,  and \textsl{parabolic} if~$\lambda=0$.

The continuous Denjoy\,--\,Wolff theorem opens the door to the study of the slopes of trajectories with respect to the Denjoy\,--\,Wolff point. Namely, given $z\in\D$, one defines the \textsl{set of slopes} $\,{\Slope^{+}(\varphi_t(z),\tau)}$ of a (hyperbolic or parabolic) semigroup $(\varphi_t)$ with associated Denjoy\,--\,Wolff point $\tau\in\partial\D$  as the $\omega$-limit of the curve
\begin{equation*}
[0,+\infty)\ni t\mapsto\mathrm{Arg(}1-\overline{\tau}\varphi_t(z))\in
\left( -\frac{\pi }{2},\frac{\pi }{2}\right) ,
\end{equation*}
where $\mathrm{Arg}$ denotes the principal argument.
It is well-known that $\Slope^{+}(\varphi_t(z),\tau)$ is a (non-empty) closed subinterval of $\left[ -\frac{\pi }{2},\frac{\pi }{2}\right]$.

If the semigroup $(\varphi_t)$ is hyperbolic, it is known \cite[Theorem 2.8]{Contreras-Diaz:pacific} that, for every $z\in\D$, $\Slope^{+}(\varphi_t(z),\tau)$ is a single point of $(-\pi/2,\pi/2)$ and
\[
\bigcup_{z\in\D}\Slope^{+}(\varphi_t(z),\tau) =(-\pi/2,\pi/2).
\]

In the parabolic case, the situation is quiet different. On the one hand, it is possible to prove that $\Slope^{+}(\varphi_t(z),\tau)$ does not depend on $z$ \cite[Theorem~2.9]{Contreras-Diaz:pacific}. On the other hand, there is no such a strong result as in the hyperbolic case and, indeed, there has been open for a long time the so-called\\[1ex]
\textbf{``The Parabolic Slope Conjecture''}:\\ \textit{
There exists $\theta\in[-\pi/2,\pi/2]$ such that $~~{\Slope^{+}(\varphi_t(z),\tau)=\{\theta\}}$ for all~${z\in\D}$.
}\vskip1ex

This conjecture has been studied and explicitly mentioned in many papers (\cite{Contreras-Diaz:pacific}, \cite{EJ}, \cite{EKRS}, \cite{ERSY}, \cite{ES}, \cite{ESY}), stating its validity in a lot of cases. To our best knowledge, the best two relevant results (of somehow different nature) are:
\begin{itemize}
\item if $(\varphi_t)$ is a parabolic one-parameter semigroup with positive hyperbolic step, then
$~\Slope^{+}(\varphi_t(z),\tau)$ is a single point and it is equal to $-\pi/2$ or $\pi/2$, see \cite{Contreras-Diaz:pacific} and~\cite{Contreras-Diaz-Pommerenke};
\item if for some $\varepsilon>0$ the limit
$$
\lim_{z\to\tau} \frac{G(z)}{(z-\tau)^{1+\varepsilon}}
$$  exists finitely and it is different from zero\footnote{Notice that for any parabolic one-parameter semigroup, $\angle\lim_{z\to\tau} \frac{G(z)}{z-\tau}
=0$.}, then $\Slope^{+}(\varphi_t(z),\tau)$ is a single point~\cite{EJ}.
\end{itemize}

In other words, with a minor assumption of regularity, the above conjecture is true. We, however, prove below that the conjecture is false in general.

\begin{theorem} There exists a parabolic semigroup $(\varphi_t)$ in the unit disk with the Denjoy\,--\,Wolff point~$\tau=1$ such that
\[
\Slope^{+}(\varphi_t(z),1)=[-\pi/2,\pi/2], \textrm{ for each } z\in\D.
\]
In particular, ``The Parabolic Slope Conjecture'' is false.
\end{theorem}

\noindent
\textbf{Note:} After finishing the redaction of this paper, Prof. D.\,Betsakos kindly communicated to the authors that he also had just obtained a counterexample to the conjecture using a different approach, based on certain estimates for the harmonic measure.

\section{Proof of the main result}

\begin{proof}
We start by introducing some notations. Given four real numbers $u_1,u_2,v,w$ with $u_1<u_2$ and $v,w>0$, we define the associated rectangle
\[
\textsc{R}(u_1,u_2,v,w):=\{z\in\C: u_1<\Re z\le u_2,-w<\Im z<v\}.
\]
This kind of rectangles will be the bricks to generate a wide family of domains which will be indexed by the set $I\times J\times J$, where:
\begin{itemize}
\item $I$ is the set of all strictly increasing sequences of positive real numbers converging to $+\infty$.
\item $J$ is the set of all increasing (not necessarily strictly) sequences of real numbers with values in $[1,+\infty)$.
\end{itemize}
Then, given $\lambda=(\alpha,\beta,\gamma)\in I\times J\times J$ with $\alpha=(u_n)$, $\beta=(v_n)$ and $\gamma=(w_n)$, we set
\[
\Omega_{\lambda}:=\textsc{R}(-1,u_1,1,1)\cup\left(\,\bigcup_{j=1}^\infty \textsc{R}(u_j,u_{j+1},v_j,w_j)\right).
\]

We note that $\Omega_{\lambda}$ is a simply connected domain of the complex plain and $[0,+\infty)\subset \Omega_{\lambda}$. This allows us to consider the canonical Riemann map $h_{\lambda}=h_{\Omega_\lambda}$ from $\D$ onto $\Omega_{\lambda}$ normalized, as usual, by
\[
h_{\lambda}(0)=0,\quad h_{\lambda}^\prime(0)>0.
\]
Moreover, by \cite[Section 9.7, Lemma 2]{Shapiro}, there exists
\[
\lim_{t\rightarrow +\infty}h_{\lambda}^{-1}(t)=:b_{\lambda}\in\partial\D.
\]
In particular, the map $\D\ni z\mapsto g_{\Omega_\lambda}(z)=g_{\lambda}(z):=h_{\lambda}(b_{\lambda}z)$ is another conformal map of~$\D$ onto $\Omega_{\lambda}$ fixing the origin. The key difference with $h_\lambda$ is that the normalization now comes from the fact
\[
\lim_{t\rightarrow +\infty}g_{\lambda}^{-1}(t)=1.
\]
Since
\[
\Omega_{\lambda}+t\subset \Omega_{\lambda}\quad \textrm{ for every } t\geq 0,
\]
we may define the mappings
\[
\varphi_t^{\Omega_{\lambda}}(z):=g_{\lambda}^{-1}(g_{\lambda}(z)+t),\quad t\geq 0,\ z\in\D.
\]
Bearing in mind the normalization of $g_{\lambda}$, we deduce that $(\varphi_t^{\Omega_{\lambda}})$ is a one-parameter semigroup of holomorphic functions in~$\UD$ having the Denjoy\,--\,Wolf point at~$1$ and with $g_{\lambda}$ as its Koenigs function.

A fundamental property of this family of domains $\mathcal{D}:=(\Omega_\lambda)_{\lambda\in I\times J\times J}$ is compiled in the following claim, whose proof will be given at the end.
\vskip1ex
\noindent{\bf Claim:} \textit{Fix $n\in\N$ and the following three vectors:
\begin{enumerate}
\item[] $(u_1,...,u_{n+1}) \textrm{ with } 0<u_1<u_2<\cdots<u_{n+1}$,
\item[] $(v_1,...,v_n) \textrm{ with } 1\leq v_1\leq v_2\leq\cdots\leq v_n$,
\item[] $(w_1,...,w_n) \textrm{ with } 1\leq w_1\leq w_2\leq\cdots\leq w_n$.
\end{enumerate}
Then, for any $\varepsilon\in(0,1)$, there exists $M_1=M_1(\varepsilon)>0$ (resp. $M_2=M_2(\varepsilon)>0$) such that for any $M\geq M_1$ (resp. $M\geq M_2$) and any $\Omega\in\mathcal{D}$ satisfying
\begin{equation}\label{EQ_Omega-cond}
\{z\in\Omega\colon \Re z\le u_{n+1}+M\}=\textsc{R}(-1,u_1,1,1)\cup\left(\,\,\bigcup_{j=1}^{n+1} \textsc{R}_j\,\,\right),
\end{equation}
where
\begin{eqnarray*}
\textsc{R}_j&:=&\textsc{R}(u_j,u_{j+1},v_j,w_j)\text{ for $j=1,\ldots,n$, and }\\
\textsc{R}_{n+1}&:=&\textsc{R}(u_{n+1},u_{n+1}+M,v_n+M,w_n)\text{ (resp. $\textsc{R}_{n+1}:=\textsc{R}(u_{n+1},u_{n+1}+M,v_n,w_n+M)$),}
\end{eqnarray*}
we have
\[
\mathrm{Arg}\left(1-g_{\Omega}^{-1}(\xi)\right)\geq \frac{\pi}{2}(1-\varepsilon)\quad\text{(resp.~~} \mathrm{Arg}\left(1-g_{\Omega}^{-1}(\xi)\right)\leq -\frac{\pi}{2}(1-\varepsilon)\text{)}
\]
for some $\xi\in(u_{n+1},u_{n+1}+M)$.}\vskip1ex

Assuming this claim, we can generate recurrently three sequences
\[
\alpha=(u_n)\in I,\ \beta=(v_n)\in J,\ \gamma=(w_n)\in J,
\]
$\lambda:=(\alpha,\beta,\gamma)$, the corresponding domain $\Omega:=\Omega_{\lambda}\in\mathcal{D}$ and an additional sequence of real numbers $(\xi_n)$ in such a way that
\begin{enumerate}
\item[(A)] $\lim_n v_n=\lim_n w_n=+\infty$,
\item[(B)] $u_{n}<\xi_n<u_{n+1}$, for every $n\in\N$,
\item[(C)] $\mathrm{Arg}\left(1-g_{\Omega}^{-1}(\xi_{2n})\right)\leq -\frac{\pi}{2}(1-\frac{1}{2n})$ for every $n\in\N$,
\item[(D)] $\mathrm{Arg}\left(1-g_{\Omega}^{-1}(\xi_{2n+1})\right)\geq \frac{\pi}{2}(1-\frac{1}{2n})$ for every $n\in\N$.
\end{enumerate}
In particular, the $\omega$-limit of the curve
\[
[0,+\infty)\ni t\mapsto\mathrm{Arg}\left(1-\varphi_t^{\Omega}(0)\right)=\mathrm{Arg}\left(1-g_{\Omega}^{-1}(t)\right)
\]
contains the numbers $\pi/2$ and $-\pi/2$. Therefore, $\Slope^{+}(\varphi^\Omega_t(0),1)=[-\pi/2,\pi/2]$. Finally, because of~(A) and appealing to
\cite[Theorem 2.9 (1)]{Contreras-Diaz:pacific}, we obtain that
\[
\Slope^{+}(\varphi^\Omega_t(z),1)=[-\pi/2,\pi/2]\quad \textrm{ for every } z\in\D.
\]
The fact that the semigroup $(\varphi^\Omega_t)$ is parabolic follows from~(A) by~\cite[Corollary 2.2]{Contreras-Diaz:pacific}.

\vspace{0.5truecm}

\noindent
\textbf{Proof of the Claim}. We only deal with one of the two possibilities (the other is completely similar) and proceed by reductio ad absurdum. So we suppose on the contrary to our claim that there exist $\varepsilon>0$, a sequence of positive real numbers $(M_k)$ with $\lim_k M_k=+\infty$ and a family of domains $(\Omega_k)$ in $\mathcal{D}$ such that, for every $k\in\N$,
\begin{enumerate}
\item[(1)] equality~\eqref{EQ_Omega-cond} holds for $\Omega:=\Omega_k$ and  $\textsc{R}_{n+1}:=\textsc{R}(u_n,u_{n+1}+M_k,v_n+M_k,w_n)$,
\item[(2)] $\mathrm{Arg}\left(1-g_{\Omega_k}^{-1}(\xi)\right)<\frac{\pi}{2}(1-\varepsilon)$ for all $\xi\in(u_{n+},u_{n+1}+M_k)$.
\end{enumerate}
It is easy to see that the sequence $(\Omega_k)$ converges to its kernel (w.r.t. the origin) in the sense of Carath{\'e}odory,
\[
\Omega_k\to\widehat{\Omega}:=\mathrm{Ker}_0 \{\Omega_k\}=\left\{z\in\C:\Re z>u_{n+1},\ \Im z>-w_n \right\}\cup\left(\,\,\bigcup_{j=1}^{n} \textsc{R}_j\,\,\right).
\]
Let $h_k$ and $h$ stand for the canonical Riemann maps of~$\UD$ onto $\Omega_k$ and $\Omega$, respectively, normalized, as usual, by $h_k(0)=h(0)=0$, $h_k'(0)>0$, $h'(0)>0$.
Then by the Carath{\'e}odory kernel convergence theorem, see, \textit{e.g.}, \cite[Chapter 15, Corollary 4.11 and Theorem 4.7]{Conway2},

\begin{enumerate}
\item[(i)] $h_k\to h$ uniformly on compacta as $k\to+\infty$,
\item[(ii)] for each compact set $K\subset \widehat\Omega$, all but a finite number of $h_k^{-1}$'s are defined on~$K$ and  $h^{-1}_k\to h^{-1}$ uniformly on~$K$ as $k\to+\infty$.
\end{enumerate}

Using again \cite[Section 9.7, Lemma 2]{Shapiro}, we obtain
\[
b_k=\lim_{t\rightarrow+\infty} h_k^{-1}(t)\in\partial\D,\ k\in\N,\quad b=\lim_{t\rightarrow+\infty} h^{-1}(t)\in\partial\D.
\]
Moreover, $b=\lim_k b_k$, essentially due to Julia's Lemma. Indeed, since $h(\D)+t\subset h(\D)$ and $h(0)=0$, we can consider the associated one-parameter semigroup of the disk ${\phi_t:=h^{-1}\circ(h+t)}$. Using (i) and (ii), we find that $\phi_t$ is the pointwise limit of~$\phi^k_t$, where $(\phi^k_t)$ is the one-parameter semigroup associated  with the couple $(\Omega_k,h_k)$, \textit{i.e.}, ${\phi^k_t:=h_k^{-1}\circ(h_k+t)}$. Now passing to the limit in Julia's inequality~\cite[(1.2.9) on p.\,49]{Abate-book}, we see that if $\mu$ is an accumulation point of $(b_k)$, then
\[
\frac{|\phi_t(z)-\mu|^2}{1-|\phi_t(z)|^2}\leq\frac{|z-\mu|^2}{1-|z|^2}, \textrm{ for all } z\in\D \textrm{ and }t\geq0,
\]
from which it follows that $\mu=b$, because $\phi_t(z)\to b$ as $t\to+\infty$ for all $z\in\UD$ according to the continuous Denjoy\,--\,Wolff theorem.

Therefore, defining $g(z):=h(bz)$, $z\in\D$, and denoting $g_k:=g_{\Omega_k}$, we also have that $g_k^{-1}\to g^{-1}$, as $k\to+\infty$,  uniformly  on compacts subsets of $\widehat{\Omega}$ in the same sense as in~(ii).

Clearly $g$ is the Koenigs function of the one-parameter semigroup~$(\psi_t)$, given by ${\psi_t(z):=\overline b\phi_t(bz)}$ for all $z\in\UD$ and $t\ge0$. By \cite[Corollary 2.2 and Theorem 2.9 (2)]{Contreras-Diaz:pacific}, we deduce that the limit
\begin{equation}\label{EQ_lim}
\theta_0:=\lim_{t\rightarrow+\infty}\mathrm{Arg}(1-g^{-1}(t))
\end{equation}
exists and must value $\pi/2$ or $-\pi/2$. We will prove now that both cases are impossible.

Assume first that $\theta_0=\pi/2$. In this case there would exist $T>u_{n+1}$ satisfying ${\mathrm{Arg}(1-g^{-1}(T))}>{\frac{\pi}{2}\left(1-\frac{\varepsilon}{2}\right)}$. Since $g^{-1}_k\to g^{-1}$ and $M_k\to+\infty$ as $k\to+\infty$, we could also find $N\in\Natural$ such that $u_{n+1}+M_N>T$ and $\mathrm{Arg}(1-g_N^{-1}(T))>\frac{\pi}{2}\left(1-{\varepsilon}\right)$, which would contradict~(2). This shows that  $\theta_0\neq\pi/2.$

Suppose now that the limit~\eqref{EQ_lim} equals $-\pi/2$.
Set $\omega:=\frac{1+i}{\sqrt{2}}$ and consider the unique conformal map~$f$ of $\UD$ onto
$
\Omega_0:=\{z\in\C: \Re z>u_{n+1}, \Im z>-w_n\}
$
such that ${f(0)=u_{n+1}-iw_n+\omega}$ and $\lim_ {u_{n+1}<t\rightarrow+\infty} f^{-1}(t)=1$.
In fact, $f$ is given by the formula
\[
f(z)=\frac{1+i}{\sqrt{2}}\sqrt{\frac{1+z}{1-z}}+u_{n+1}-iw_n,\ z\in\D,
\]
where the square root is to be understood as its principle branch. In particular, it is easy to check that
\begin{equation}\label{EQ_slope-f}
\lim_{u_{n+1}<t\rightarrow+\infty}\mathrm{Arg}(1-f^{-1}(t))=\frac{\pi}{2}.
\end{equation}

By construction, $f(\D)\subset\widehat{\Omega}$. Moreover, $\widehat{\Omega}$ and $f(\D)$ are Jordan domains, so $g$ and $f$ both extend homeomorphically to the closed unit disk. Therefore, $\phi:=g^{-1}\circ f$ is a holomorphic self-map of~$\UD$ that also extends homeomorphically to~$\overline{\UD}$.  Recall that $\lim_{t\rightarrow+\infty}g^{-1}(t)=1$. Hence $\phi(1)=1$. Moreover, bearing in mind the shape of the (Jordan) domains $f(\D)$ and~$\widehat{\Omega}$, we can apply the Schwarz's reflection principle to extend $\phi$ holomorphically through an open arc of~$\partial\UD$ containing the point $1$. Note that the extended map is univalent. It follows that
\begin{equation}\label{EQ_lim-phi}
\frac{1-\phi(z)}{1-z}\to\phi'(1)>0\text{~~~as~~$\UD\ni z\to 1$}\quad(\textrm{the unrestricted limit in }\D).
\end{equation}
Finally, consider the following equality, for $t>0$ large enough,
\[
\frac{1-g^{-1}(t)}{|1-g^{-1}(t)|}=\frac{1-\phi(f^{-1}(t))}{1-f^{-1}(t)}\frac{1-f^{-1}(t)}{|1-f^{-1}(t)|}
\frac{|1-f^{-1}(t)|}{|1-\phi(f^{-1}(t))|}
\]
On the one hand, since the limit~\eqref{EQ_lim} is $-\pi/2$, the left-hand side must tend to~$-i$ as $t\to+\infty$. On the other hand, by \eqref{EQ_slope-f} and~\eqref{EQ_lim-phi},  the right-hand side tends to~$i$ as~$t\to+\infty$. This contradiction completes the proof of the Claim.
\end{proof}


\begin{thebibliography}{99}
\bibitem{Abate-book} M.\,Abate, Iteration theory of holomorphic maps on taut manifolds.
Research and Lecture Notes in Mathematics. Complex Analysis and Geometry,
Mediterranean, Rende, 1989. MR1098711 (92i:32032).

\bibitem{ABCD} M.\,Abate, F.\,Bracci, M.D.\,Contreras, and S.\,D\'iaz-Madrigal, The evolution of Loewner's differential
equations, Newsletter European Math. Soc. {\bf 78} (December) (2010), 31--38.

\bibitem{BP} E.\,Berkson\ and\ H.\,Porta, Semigroups of analytic functions and composition
operators. Michigan Math. J. {\bf 25} (1978), no.~1, 101--115. MR0480965 (58
\#1112)

\bibitem{Contreras-Diaz:pacific} M.\,D.\,Contreras and S.\,D\'{\i}az-Madrigal, Analytic flows in the unit disk: angular derivatives and boundary fixed points.
Pacific J. Math. \textbf{222} (2005), 253--286.


\bibitem{Contreras-Diaz-Pommerenke} M.\,D.\,Contreras, S.\,D\'{\i}az-Madrigal, and Ch.\,Pommerenke,  Some Remarks on the Abel equation in the unit disk, J. London Math. Soc. (2), \textbf{75}, (2007), 623--634.

\bibitem{Conway2} J.\,B.\,Conway, Functions of One Complex Variable,~II.
Second edition. Graduate Texts in Mathematics, 159. Springer-Verlag, New York-Berlin, 1996.

\bibitem{EJ} M.\,Elin and F.\,Jacobson, Parabolic type semigroups: asymptotics and order of contact. Preprint 2013. Available at  http://arxiv.org/pdf/1309.4002.pdf.

\bibitem{EKRS} M.\,Elin, D.\,Khavinson, S.\,Reich, and D.\,Shoikhet, Linearization models for parabolic dynamical systems via Abel's functional equation, Ann. Acad. Sci. Fen. {\bf 35} (2010), 1--34.

\bibitem{ERSY} M.\,Elin, S.\,Reich, D.\,Shoikhet, and F.\,Yacobson, Asymptotic behavior of one-parameter semigroups and rigidity of holomorphic generators, Complex Anal. Oper. Theory {\bf 2} (2008), 55--86.

\bibitem{ES-book} M.\,Elin and D.\,Shoikhet, Linearization models for complex dynamical systems. Topics in univalent functions, functional equations and semigroup theory. Birkh\"auser Basel, 2010.

\bibitem{ES} M.\,Elin and D.\,Shoikhet, Boundary behavior and rigidity of semigroups of holomorphic mappings, Anal. Math. Phys. {\bf 1} (2011), 241--258.

\bibitem{ESY} M.\,Elin,  D.\,Shoikhet, and F.\,Yacobson, Linearization models for parabolic type semigroups, J. Nonlinear Convex Anal.  {\bf 9} (2008), 205--214.

\bibitem{Goryainov::survey}V.\,V.\,Gorya\u\i nov, Semigroups of analytic functions in analysis and applications. Uspekhi Mat. Nauk {\bf 67} (2012), no. 6(408), 5--52; translation in Russian Math. Surveys {\bf 67} (2012), no.~6, 975--1021. MR3075076

\bibitem{Poggi-Corradini}P.\,Poggi-Corradini, Backward-iteration sequences with bounded hyperbolic steps for analytic self-maps of the disk. Rev. Mat. Iberoamericana {\bf 19} (2003), 943--970.

\bibitem{Pommerenke} Ch.\,Pommerenke, On the iteration of analytic
functions in a half plane I, J. London Math. Soc. (2) \textbf{19} (1979),
439--447.

\bibitem{Shapiro} J.H.\,Shapiro, Composition Operators and Classical
Function Theory. Springer-Verlag, New York, 1993.

\bibitem{Shoikhet-book} D.\,Shoikhet, Semigroups in Geometrical Function Theory. Kluwer Academic Publishers, Dordrecht, 2001.

\end{thebibliography}
\end{document}